\documentclass[a4paper,11pt]{article}
\usepackage{mathrsfs}
\usepackage{latexsym}
\usepackage{amsmath,amssymb}
\usepackage{amsthm}
\usepackage{amsfonts}
\usepackage[usenames]{color}
\usepackage{amssymb}
\usepackage{graphicx}
\usepackage{amsmath}
\usepackage{amsfonts}
\usepackage{amsthm}
\usepackage{mathrsfs}
\usepackage{dsfont}
\usepackage{indentfirst}

\newcommand\pa{\partial}

\newcommand\hlf{\frac{1}{2}}

\newcommand\Si{\Sigma}

\def\p{\partial }
\def\eps{\varepsilon}

\def\phi{\varphi}
\def\lam{\lambda}
\def\ga{\gamma}

\newtheoremstyle{mythm}{1.5ex plus 1ex minus .2ex}{1.5ex plus 1ex
minus .2ex}{\kai}{\parindent}{\song\bfseries}{}{1em}{}
\numberwithin{equation}{section}

\newtheorem{theorem}{Theorem}[section]
\newtheorem{lemma}{Lemma}[section]

\newtheorem{proposition}{Proposition}[section]
\newtheorem{remark}{Remark} [section]

\allowdisplaybreaks[4]

\begin{document}
\title
{Liouville type theorem of integral equation with anisotropic structure}
\author{Yating Niu}
\date{}
\maketitle

\begin{abstract}
In this paper, we classify all positive solutions for the following integral
equation:
\begin{equation}\label{s0}
u(x)=\int_{\mathbb{R}^n_+}K_b(x,y)y_n^b f(u(y))dy,
\end{equation}
where $ b > 1$ is a constant. Here $ K_b(x,y)$ is the Green function of the following homogeneous Neumann boundary problem
\begin{equation}\label{s11}
\left\{
\begin{aligned}
-\text{div}(x^{b}_n \nabla u)&= f \quad \text{in} \ \mathbb R^n_+\\
\frac{\p u}{\p x_n}&= 0 \quad \text{on} \ \p \mathbb R^n_+ .
\end{aligned}
\right.
\end{equation}
By using the method of moving planes in integral form, we derive the symmetry of positive solutions. We also establish the equivalence between the integral equation and its corresponding partial differential equation. Similarly, the results can be generalized to the integral system.
\end{abstract}

\section{Introduction}
In this paper, we study the positive solutions $ u(x)$ of the following type of integral equation
\begin{equation}\label{s3}
u(x)=\int_{\mathbb R^n_+}K_b(x,y)y_n^b f(u(y))dy
\end{equation}
and integral system
\begin{equation}\label{p2}
\left\{
\begin{aligned}
u(x)&=\int_{\mathbb R^n_+}K_b(x,y)y_n^b f(v(y))dy \quad x\in \mathbb R^n_+\\
v(x)&=\int_{\mathbb R^n_+}K_b(x,y)y_n^b g(u(y))dy \quad x\in \mathbb R^n_+.
\end{aligned}
\right.
\end{equation}
on upper half space $ \mathbb R^n_+$, where $ n \geq 3$.

In recent years, there has been great interest in using the method of moving planes to classify the solutions of equation. It is a very powerful tool to study the symmetry of solution. The method of moving planes for PDEs was invented by the Alexandroff in the early 1950's. Later, it was further developed by Serrin \cite{Serrin} and Gidas, Ni and Nirenberg \cite{GiNiNi1979}, \cite{GidasNiNiren}. In the paper of Chen, Li \cite{ChenLi91}, they studied the following partial differential equation
\begin{equation}\label{k1}
\Delta u + u^p =0 \quad x\in \mathbb R^n.
\end{equation}
They proved that for $ p=\frac{n+2}{n-2}$ all the positive $ C^2$ solutions of \eqref{k1} are radially symmetric about some point.

It is a natural question to ask whether similar results hold for the following system
\begin{equation}\label{k2}
\left\{
\begin{aligned}
-\Delta u &= v^p  \quad x\in \mathbb{R}^n \\
-\Delta v &= u^q  \quad x\in \mathbb{R}^n.
\end{aligned}
\right.
\end{equation}
In 1993, Mitidieri \cite{EM93} has considered system \eqref{k2} the nonexistence of radial positive solutions, where $ \frac{1}{p+1} + \frac{1}{q+1} > \frac{n-2}{n}$ and $ p >1$, $ q >1$. Later, Mitidieri in \cite{EM} also extended this results to more general system. In \cite{DeFe}, they proved that for $ p=q=\frac{n+2}{n-2}$ the positive solutions of \eqref{k2} are radially symmetric with respect to some point. Guo and Liu extended to more general elliptic system in \cite{GuoLiu}. For other results, we refer to \cite{PRE,GuoNie,Huang2012,HuangLi,SeZou}.

For the integral equation, we can use the method of moving planes in integral form to study the properties of the solution. The integral equation \eqref{s3} and integral system \eqref{p2} is closely related to \cite{Yu13}. In this paper, Yu studied the positive solutions for the following integral equation
\begin{equation}\label{s2}
u(x) = \int_{\mathbb R^n} \frac{1}{|x-y|^{n-\alpha}} f(u(y)) dy  \quad x\in \mathbb R^n,
\end{equation}
where $ n \geq 2$, $ 0<\alpha<n$. If $ f(u)= u ^{\frac{n+\alpha}{n-\alpha}}$, integral equation \eqref{s2} arises as an Euler-Lagrange equation for a functional under a constraint in the context of Hardy-Littlewood-Sobolev inequalities. Lieb \cite{Lieb} posed the classification of all the critical points of the functional - the solutions of the integral equation as an open problem. Chen, Li and Ou \cite{ChenLiOu06} solved the problem by using the method of moving planes in an integral form. They proved that all the positive regular solutions are radially symmetric and monotone decreasing about some point.

For the integral system, Yu also established the Liouville type result for the following integral system
\begin{equation}\label{p3}
\left\{
\begin{aligned}
u(x)&=\int_{\mathbb R^n}\frac{1}{|x-y|^{n-\alpha}} f(u(y),v(y))dy  \quad x\in \mathbb R^n\\
v(x)&=\int_{\mathbb R^n}\frac{1}{|x-y|^{n-\alpha}} g(u(y),v(y))dy  \quad x\in \mathbb R^n
\end{aligned}
\right.
\end{equation}
in \cite{Yu13}, where $ n \geq 2$, $ 0<\alpha<n$. If $ f=v^q$ and $ g=u^p$, Chen, Li and Ou \cite{ChenLiOu05} proved that the positive solutions of \eqref{p3} are radially symmetric for $ \frac{1}{q+1} + \frac{1}{p+1}=\frac{n-\alpha}{n}$. Later, more general integral equations and systems have also been studied in papers \cite{ChenLiOu09,GuoLiu,GuoNie,LiZhuo,LiuGuoZhang,MaChen,EM}.

The above results are all about the whole space $ \mathbb{R}^n$. For the upper half space
\[
\mathbb{R}^n_+ = \{x=(x_1,x_2,\cdots,x_n)\in \mathbb{R}^n \ | \ x_n>0\},
\]
Fang and Chen \cite{FangChen} considered the integral equation $ u(x)= \int_{\mathbb R^n_+} G^+_\infty(x,y)u^p dy$, where $ G^+_\infty(x,y)$ is the Green's function in $ \mathbb R^n_+$ with the following Dirichlet problem
\begin{equation}\label{s12}
\left\{
\begin{aligned}
(-\Delta)^m u &= u^p \qquad \qquad \qquad \quad \ \text{in} \ \mathbb R^n_+\\
u=\frac{\p u}{\p x_n}&=\cdots = \frac{\p^{m-1}u}{\p x^{m-1}_n} =0  \quad   \text{on} \ \p\mathbb R^n_+.
\end{aligned}
\right.
\end{equation}
They proved that the Dirichlet problem \eqref{s12} is equivalent to the integral equation. Later, Tang and Dou \cite{TangDou} studied the system of integral equations on upper half space. In 2015, Chen, Fang and Yang \cite{ChenFangYang} considered the Dirichlet problem involving the fractional Laplacian on upper half space.

In this paper, we study the upper half space results for problem \eqref{s3} and \eqref{p2}. Our result is the following
\begin{theorem}\label{thm0}
Let $ u(x)\in C(\overline{\mathbb{R}^n_+})$ be a positive solution of \eqref{s3}, where $ n\geq3$ and let $ f:[0,+\infty)\rightarrow \mathbb{R}$ is a continuous function with the properties\\
$(i) \ f(t)$ is non-decreasing in $ [0,+\infty)$; \\
$(ii)g(t)=\frac{f(t)}{t^{\frac{n+b+2}{n+b-2}}}$ is non-increasing in $ (0,+\infty)$.\\
Then either $ u\equiv 0$ or there exists $ x_0 \in \pa\mathbb{R}^{n}_+$ such that $ u(x)=u_{a,x_0}(x)=\left(\frac{ca}{a^2+|x-x_0|^2}\right)^{\frac{n+b-2}{2}}$ and $ g(t)\equiv \bar{c}$, where $ c=\sqrt{\frac{(n+b)(n+b-2)}{\bar{c}}}$ and $ a\geq0 $.
\end{theorem}
\begin{remark}\label{r1}
In Theorem $ \ref{thm0}$, we can get $ f(t)\geq 0$ from conditions $ (i)$ and $ (ii)$.
\end{remark}

\begin{theorem}\label{thm1}
Let $ (u,v) \in C(\overline{\mathbb R^n_+})\times C(\overline{\mathbb R^n_+}) $ be a positive solution of problem \eqref{p2}. Suppose that $ f$, $ g$: $ [0,+\infty) \rightarrow \mathbb R$ are continuous and satisfy \\
(i) \ $ f(t)$, $ g(t)$ are non-decreasing in $ [0,+\infty)$; \\
(ii) \ $ h(t)=\frac{f(t)}{t^{\frac{n+b+2}{n+b-2}}}$, $ k(t)=\frac{g(t)}{t^{\frac{n+b+2}{n+b-2}}}$ are non-increasing in $ (0,+\infty)$. \\
Then either $ (u,v)\equiv(0,0) $ or $ u$, $ v$ has the form $ u(x)=\left(\frac{ca}{a^2+|x-x_0|^2}\right)^{\frac{n+b-2}{2}}$ and $ v(x)=\left(\frac{\tilde{c}a}{a^2+|x-x_0|^2}\right)^{\frac{n+b-2}{2}}$, for some $ x_0 \in \pa\mathbb{R}^{n}_+$ and $ a\geq 0$. \end{theorem}

\begin{remark}
In Theorem $ \ref{thm1}$, we can get $ f(t)\geq 0$ and $ g(t)\geq 0$ from conditions $ (i)$ and $ (ii)$.
\end{remark}

The main method in this paper is the moving plane method in integral form. We also use Hardy-Littlewood-Sobolev inequalities and Kelvin's transform to prove the symmetry of positive solutions with respect to the $ x_1, \cdots ,x_{n-1}$ directions. And $ x_n$ is the anisotropic direction.

This paper is organized as follow. In Section 2, we use the method of moving planes to get the symmetry of the solutions. In Section 3, we prove the equivalence of the integral equation with differential equation and the nonexistence of an ordinary differential equation. We also give the proof of Theorem \ref{thm0}. In Section 4, we obtain the symmetry of solutions of integral equation system and and we prove Theorem $ \ref{thm1}$.

\section{ Symmetry of solutions}

Let us study the positive solutions to the integral equation
\[
u(x)=\int_{\mathbb R^n_+}K_b(x,y)y_n^b f(u(y))dy \quad x \in \mathbb R^n_+ ,
\]
where
\begin{equation}\label{G1}
\begin{split}
&K_b(x,y)=D_b \int_0^1 (|x-y|^2(1-\tau)+|x-y^*|^2\tau)^{-\frac{n+b-2}{2}}[\tau(1-\tau)]^{\frac b 2-1}d\tau, \\
&D_b=2^{b-2}\pi^{-\frac n2}\frac{\Gamma(\frac{n+b-2}{2})}{\Gamma(\frac b2)},\quad y^*=(y_1,\cdots,y_{n-1},-y_n).
\end{split}
\end{equation}
Moreover, for  $b>0,n\geq 3$, the following estimates for $K_{b}(x,y)$ holds
\begin{equation}\label{G2}
|\p_x^{\gamma} K_b(x,y)|\leq C(\gamma,b) |x-y^*|^{-b}|x-y|^{2-n-|\gamma|},
\end{equation}
where $\gamma \in \mathbb N^n $ (see Proposition 1-2 of \cite{Horiuchi95}).
Since we don't know the behaviors of $ u$ at infinity, we introduce the Kelvin's transform of $ u$ as $ v(x)=\frac{1}{|x-x_0|^{n+b-2}}u\left(\frac{x-x_0}{|x-x_0|^2}\right)$.

$ \forall x_0 \in \p\mathbb{R}^n_+ $, we define $ w(x):= u(x-x_0)$. Since
\begin{align*}
w(x) &= u(x-x_0) \\
&= \int_{\mathbb{R}^n_+} K_b(x,y+x_0)y^b_n f(u(y)) dy \\
&= \int_{\mathbb{R}^n_+} K_b(x,y)y^b_n f(w(y)) dy,
\end{align*}
so we take $ x_0 = 0$. We consider $ v(x)=\frac{1}{|x|^{n+b-2}}u\left(\frac{x}{|x|^2}\right)$. Then a direct calculation shows that $ v(x)$ solves \\
\[
v(x)=\int_{\mathbb R^n_+}K_b(x,y)y_n^b f(v(y)|y|^{n+b-2})\frac{1}{|y|^{n+b+2}}dy \quad x \in \mathbb R^n_+ .
\]
We define $ \tilde{\tau}=\frac{n+b+2}{n+b-2}$. By the definition of \(g(t)=\frac{f(t)}{t^{\tilde{\tau}}}\displaystyle\), we deduce that $ v(x)$ satisfies \\
\[
v(x)=\int_{\mathbb R^n_+}K_b(x,y)y_n^b g(v(y)|y|^{n+b-2}) v^{\tilde{\tau}}dy .
\]

Since $ u$ is continuous in $ \overline{\mathbb{R}^n_+}$, we conclude that $ v$ is continuous and positive in $ \overline{\mathbb{R}^n_+} \setminus \{0\}$ with possible singularity at the origin. Moreover, $ v$ decays at infinity as $ u(0)|x|^{2-n-b}$. By the asymptotic behavior of $ v$ at $ \infty$, we get
\[
v(x) \in L^{\frac{2n}{n+b-2}}(\mathbb{R}^n_+ \setminus B_r(0)) \cap L^\infty (\mathbb{R}^n_+ \setminus B_r(0))
\]
for all $ r>0$. Now we use the moving plane method to prove our result.

For a given real number $ \lambda$, define
\[
\Si_\lam = \{x=(x_1,\cdots ,x_n)\in \mathbb R^n_+ \ | \ x_1\geq \lambda \}, \quad T_\lambda = \{x\in \mathbb R^n_+ \ | \ x_1=\lambda \}
\]
and let $ x^\lambda = (2\lambda-x_1,x_2,\cdots ,x_n)$ and $ u_\lambda(x)=u(x^\lambda)$.
\begin{lemma}\label{lem1}
\begin{align*}
v(x)-v(x^\lambda)=& \int_{\Sigma_\lambda}(K_b(x,y)-K_b(x^\lambda,y)) y^b_n \times \\ &[g(v(y)|y|^{n+b-2})v(y)^{\tilde{\tau}}-g(v(y^\lambda)|y^\lambda|^{n+b-2})v(y^\lambda)^{\tilde{\tau}}]dy.
\end{align*}
\end{lemma}
\noindent \emph{Proof.} By the definition of $ K_b(x,y)$, we know $ K_b(x,y) = K_b(x^{\lam},y^{\lam})$, $ K_b(x^{\lam},y) = K_b(x,y^{\lam}) $. A direct calculation implies
\begin{equation*}
\begin{aligned}
v(x) =  \int_{\Sigma_\lambda} K_b(x,y) y^b_n g(v(y)|y|^{n+b-2})v(y)^{\tilde{\tau}}dy
+ \int_{\Sigma_\lambda} K_b(x,y^{\lam}) y^b_n
g(v(y^\lambda)|y^\lambda|^{n+b-2})v(y^\lambda)^{\tilde{\tau}}dy,
\end{aligned}
\end{equation*}
and
\begin{align*}
v(x^\lambda)= \int_{\Sigma_\lambda} K_b(x^{\lam},y) y^b_n
g(v(y)|y|^{n+b-2})v(y)^{\tilde{\tau}} dy
+  \int_{\Sigma_\lambda} K_b(x^{\lam},y^{\lam})
y^b_n g(v(y^\lambda)|y^\lambda|^{n+b-2}) v(y^\lambda)^{\tilde{\tau}}dy.
\end{align*}
We get the desired result. \hfill $ \Box$

\begin{lemma}\label{lem2}
Under the assumptions of Theorem \ref{thm0}, there exists $ \lam_0 > 0$ such that for all $ \lam \geq \lam_0$, we have $ v_\lam(x) \geq v(x)$ for all $ x \in \Si_\lam$.
\end{lemma}
\noindent \emph{Proof.} \ $ \forall \lam >0$, it is easy to see that $ |y|>|y^\lam|$, $ \forall y\in \Si_\lam$. Let
\[
\Si^{-}_{\lam} = \{y\in \Si_\lam \ | \ v(y)>v_\lam(y) \},
\]
then $ \forall y\in \Si^-_\lam$, we have
\[
v(y)|y|^{n+b-2} > v_\lam(y)|y^\lam|^{n+b-2}.
\]
Since $ g(t)$ is non-increasing (see Theorem \ref{thm0}(ii)), we get
\[
g(v(y)|y|^{n+b-2}) \leq g(v_\lam(y)|y^\lam|^{n+b-2}).
\]
This implies
\begin{align}
g(v(y)|y|^{n+b-2})v(y)^{\tilde{\tau}}- g(v_\lambda(y)|y^\lambda|^{n+b-2})v_\lambda(y)^{\tilde{\tau}}  \leq g(v(y)|y|^{n+b-2})(v(y)^{\tilde{\tau}} - v_\lam(y)^{\tilde{\tau}}). \label{s4}
\end{align}
As for $ y \in \Si_\lam\setminus\Si^-_\lam$, using the assumptions of Theorem \ref{thm0} we have
\begin{align}
g(v(y)|y|^{n+b-2})v(y)^{\tilde{\tau}} &= \frac{f(v(y)|y|^{n+b-2})}{|y|^{n+b+2}} \nonumber \\
&\leq \frac{f(v_\lam(y)|y|^{n+b-2})}{|y|^{n+b+2}} \nonumber  \\
&=\frac{f(v_\lam(y)|y|^{n+b-2})}{[|y|^{n+b-2}v_\lam(y)]^{\tilde{\tau}}}v_\lambda(y)^{\tilde{\tau}}\nonumber \\
&= g(v_\lam(y)|y|^{n+b-2}) v_\lam(y)^{\tilde{\tau}} \nonumber \\
&\leq g(v_\lam(y)|y^\lam|^{n+b-2})v_\lam(y)^{\tilde{\tau}}.\label{s5}
\end{align}
Then by Lemma \ref{lem1} and the mean value theorem, we know that there exists a $ \xi$ between $y$ and $ y^\lambda$ that makes
\begin{equation}
v(y)^{\tilde{\tau}}-v_\lambda(y)^{\tilde{\tau}}=\tilde{\tau}v(\xi)^{\tilde{\tau}-1}(v(y)-v_\lambda(y))
\end{equation}
true. Combining Lemma \ref{lem1}, \eqref{s4} and \eqref{s5}, we have
\begin{align*}
v(x)-v_\lambda(x)
&\leq \int_{\Sigma^-_\lambda} (K_b(x,y)-K_b(x^\lambda,y))y^b_n g(v(y)|y|^{n+b-2})(v(y)^{\tilde{\tau}}-v_\lambda(y)^{\tilde{\tau}}) dy \\
&\leq \tilde{\tau} \int_{\Sigma^-_\lambda} (K_b(x,y)-K_b(x^\lambda,y))y^b_n g(v(y)|y|^{n+b-2})v(y)^{\tilde{\tau}-1}(v(y)-v_\lambda(y)) dy.
\end{align*}

In getting the above inequality, we also used $ \forall x,y \in \Sigma_\lam$, $ K_b(x,y) \geq K_b(x^\lam, y)$ and $ \forall y \in \Sigma^-_\lam$, $ v_\lam (y) \leq v(\xi) \leq v(y)$. Since
\[
|y|^{n+b-2}v(y)=u\left(\frac{y}{|y|^2}\right) \geq \min_{B_{\frac{1}{\lam}}} u(x) \geq C_\lam >0, \quad \forall \lam>0, \ y \in \Sigma_\lam,
\]
and $ f$ is continuous, we conclude that $ g(|y|^{n+b-2}v(y))$ is bounded for $ y\in \Sigma_\lambda$. Hence we can deduce from the above inequality that
\[
v(x)-v_\lambda(x) \leq C\int_{\Sigma^-_\lambda} K_b(x,y)y^b_n v(y)^{\tilde{\tau}-1}(v(y)-v_\lambda(y)) dy, \quad \forall \lam > 0, \ x\in \Sigma_\lam.
\]
By \eqref{G2} and since $ \frac{y^b_n}{|x-y*|^b} \leq 1$, then
\[
v(x)-v_\lambda(x) \leq C\int_{\Sigma^-_\lambda} \frac{1}{|x-y|^{n-2}} v(y)^{\tilde{\tau}-1}(v(y)-v_\lambda(y)) dy.
\]
By the Hardy-Littlewood-Sobolev inequality, it follows that for any $ q > \frac{n}{n-2}$,
\[
\|v-v_\lam \|_{L^q(\Sigma^-_\lambda)} \leq C \|v^{\tilde{\tau} -1}(v-v_\lam )\|_{L^{\frac{qn}{n+2q}}(\Sigma^-_\lambda)}.
\]
By using the generalized H$ \ddot{o}$lder inequality, we obtain
\begin{align*}
\|v-v_\lam \|_{L^q(\Sigma^-_\lambda)}
& \leq C\|v^{\tilde{\tau} -1}\|_{L^{\frac{n}{2}}(\Sigma^-_\lambda)} \|v-v_\lam \|_{L^q(\Sigma^-_\lambda)}\\
&=C\left(\int_{\Sigma^-_\lambda}v^{\frac{2n}{n+b-2}}(y)dy\right)^{\frac{2}{n}}\|v-v_\lambda \|_{L^q(\Sigma^-_\lambda)}.
\end{align*}
Due to $ v\in L^{\frac{2n}{n+b-2}}(\mathbb{R}^n_+\setminus B_r(0))$, we can choose $ \lambda_0$ sufficiently large, such that for $ \lambda \geq \lambda_0$, we have
\[
C\left(\int_{\Sigma_\lambda}v^{\frac{2n}{n+b-2}}(y)dy \right)^{\frac{2}{n}} \leq \frac{1}{2}.
\]
Then we conclude
\[
\|v-v_\lam \|_{L^q(\Sigma^-_\lambda)} = 0
\]
for all  $ \lambda \geq \lambda_0$. Thus, $ \Si^-_\lam$ must be measure $ 0$. Since $ v$ is continuous, we deduce that  $ \Si^-_\lam$ is empty. \hfill $ \Box$

We define $ \lam_1 = \inf \{ \lam \ | \ v(x)\leq v_\mu(x), \  \forall \mu \geq \lam, \  x\in\Sigma_\mu \}$.
\begin{lemma}\label{lem3}
If $ \lambda_1 >0$, then $ v(x)\equiv v_{\lambda_1}(x)$ for all $ x\in \Sigma_{\lambda_1}$.
\end{lemma}
\noindent \emph{Proof.} \ Suppose that conclusion does not hold. We have $ v(x) \leq v_{\lam_1}(x)$, but $ v(x) \not\equiv v_{\lam_1}(x)$ in $ \Si_{\lam_1}$. We will prove that the plane can be moved further to the left. We infer from the proof of Lemma \ref{lem2} that
\begin{equation}\label{s6}
\|v-v_\lam \|_{L^q(\Si^-_\lam)} \leq C\left(\int_{\Si^-_\lam}v^{\frac{2n}{n+b-2}}(y)dy\right)^{\frac{2}{n}}\|v-v_\lam \|_{L^q(\Si^-_\lam)}.
\end{equation}
If one can show that for $ \eps$ sufficiently small so that $ \forall \lam \in (\lam_1-\eps,\lam_1]$, there holds
\begin{equation}\label{s7}
C\left(\int_{\Si^-_\lam}v^{\frac{2n}{n+b-2}}(y)dy\right)^{\frac{2}{n}} \leq \hlf,
\end{equation}
then by \eqref{s6}, we have $ \|v-v_\lam\|_{L^q(\Si^-_\lam)}=0$, and therefore $ \Si^-_\lam$ must be measure zero. Since $ v$ is continuous, we deduce that  $ \Si^-_\lam$ is empty. This contradicts the definition of $ \lam_1$.

Now we verify inequality \eqref{s7}. Since $ v\in L^{\frac{2n}{n+b-2}}(\mathbb{R}^n_+\setminus B_r(0))$, for any small $ \eta >0$ we can choose $ R$ sufficiently large so that
\begin{equation}\label{s8}
C\left(\int_{\mathbb{R}^n_+ \backslash B_R}v^{\frac{2n}{n+b-2}}(y)dy\right)^{\frac{2}{n}} \leq \eta.
\end{equation}
We fix this $ R$ and then show that the measure of $ \Si^-_\lam \cap B_R$ is sufficiently small for $ \lam$ close to $ \lam_1$. Since the measure of $ \Si^-_{\lam_1}$ is zero, by Lemma \ref{lem1} and \eqref{s5} we deduce $ v(x) < v_{\lam_1}(x)$ in the interior of $ \Si_{\lam_1}$.

For any $ \ga >0$, let
\[
E_\ga = \{x\in\Si_{\lam_1} \cap B_R \ | \ v_{\lam_1}(x) - v(x) > \ga\}, \qquad F_\ga = (\Si_{\lam_1} \cap B_R)\backslash E_\ga.
\]
It is obviously that
\[
\lim_{\ga\rightarrow 0} \mu(F_\ga) =0.
\]
For $ \lam < \lam_1$, let
\[
D_\lam = (\Si_\lam \backslash \Si_{\lam_1}) \cap B_R.
\]
Apparently, the measure of $ D_\lam$ is small for $ \lam$ to close to $ \lam_1$. Then it is easy to see that
\begin{equation}\label{s9}
(\Si^-_\lam \cap B_R) \subseteq (\Si^-_\lam \cap E_\ga) \cup F_\ga \cup D_\lam.
\end{equation}
In fact, $ \forall x\in \Si_\lam^- \cap E_\ga$, we have
\[
v(x)- v_\lam(x) = v(x)-v_{\lam_1}(x) + v_{\lam_1}(x) - v_\lam(x) >0.
\]
Hence
\[
v_{\lam_1}(x) - v_\lam(x) > v_{\lam_1}(x)- v(x) > \ga.
\]
It follows that
\begin{equation}\label{s10}
(\Si_\lam^- \cap E_\ga) \subseteq G_\ga \equiv \{x\in B_R \ | \ v_{\lam_1}(x) - v_\lam(x) >\ga \}.
\end{equation}
By the well-known Chebyshev inequality, we have
\begin{align*}
\mu (G_\ga) & \leq \frac{1}{\ga^{p+1}} \int_{G_\ga} |v_{\lam_1}(x)-v_\lam(x)|^{p+1}dx \\
& \leq \frac{1}{\ga^{p+1}} \int_{B_R} |v_{\lam_1}(x)-v_\lam(x)|^{p+1}dx,
\end{align*}
where $ p>0$. For each fixed $ \ga$, as $ \lam$ is close to $ \lam_1$, the right hand side of the above inequality can be made as small as we wish. Therefore by \eqref{s9} and \eqref{s10}, the measure of $ \Si^-_\lam \cap B_R$ can also be made sufficiently small. Combining this with \eqref{s8}, we obtain \eqref{s7}. \hfill $ \Box$

\section{ Proof for Theorem 1.1 }
We define
\[
\mathcal L_b(u)=- \text{div}(y_n^b \nabla u).
\]
Take a cutoff function $ \varphi_R (y) \in C^{\infty}_c(B_R)$ such that $ 0\leq \varphi_R \leq 1$ in $ B_R$, $ \varphi_R =1$ in $ B_{\frac{R}{2}}$. Setting
\[
u_R(x)=\int_{\mathbb R^n_+} K_b(x,y)y^b_n f(u(y))\varphi_R (y)dy,
\]
we clearly have $ f(u(y))\varphi_R (y) \in C(\overline{\mathbb R^n_+}) \cap \mathcal{E}'$, where $ \mathcal{E}'$ is the dual space of $ C^\infty(\overline{\mathbb R^n_+})$. By the result of \cite{Horiuchi95}, we have $ \mathcal L_b K_b(x,y)=\delta(x-y)$ and
\begin{equation}
\left\{
\begin{aligned}
\mathcal L_b(u_R) &= y^b_n f(u(y))\varphi_R (y) \quad \text{in} \quad \mathbb R^n_+ \\
\frac{\p u_R}{\p x_n} &= 0 \quad \qquad \qquad \qquad \  \text{on} \quad \p\mathbb R^n_+
\end{aligned}
\right.
\end{equation}
(see Proposition 1-1 of \cite{Horiuchi95}). By letting $ R \rightarrow \infty$, we then conclude that
\begin{equation}
\left\{
\begin{aligned}
-(y^b_n \Delta u + b y^{b-1}_n u_n)&= y^b_n f(u(y)) \quad \text{in} \quad \mathbb R^n_+ \\
\frac{\p u}{\p x_n} &= 0 \quad \qquad \quad \ \ \text{on} \quad \p\mathbb R^n_+.
\end{aligned}
\right.
\end{equation}
Therefore, the integral equation \eqref{s3} satisfies the following partial differential equation:
\begin{equation}\label{s13}
\left\{
\begin{aligned}
\Delta u + \frac{b}{y_n}u_n + f(u(y))&= 0 \quad \text{in}\quad \mathbb R^n_+ \\
\frac{\p u}{\p x_n}&=0 \quad \text{on}\quad \p \mathbb R^n_+.
\end{aligned}
\right.
\end{equation}
Then we will establish the equivalence between the integral equation \eqref{s3} and the  partial differential equation \eqref{s13}. We need the asymptotic behavior $ u(x) \sim \frac{1}{|x|^{n+b-2}}$ as $ x\rightarrow \infty$; that is, there exist two constant $ R$, $ C$ such that
\[
u(x)= \frac{C}{|x|^{n+b-2}}, \quad |x|>R.
\]
\begin{lemma}\label{lem8}
Let $ f(t)=\bar{c}t^{\frac{n+b+2}{n+b-2}}$ in Theorem \ref{thm0} and  $ u(x) \sim \frac{1}{|x|^{n+b-2}}$ as $ x\rightarrow \infty$. Then the positive solution $ u$ of $ \eqref{s3}$ is  $ C^2(\overline{\mathbb{R}^n_+})$.
\end{lemma}
\noindent $ \emph{Proof.}$  \\
\textbf{Step 1:} \emph{We prove that $ u \in C^\alpha(\overline{\mathbb{R}^n_+})$. }

Since $ u \in C^0(\overline{\mathbb{R}^n_+})$, one know $ u\in L^\infty(\overline{B_R \cap \mathbb{R}^n_+})$. Since $ u(x) \sim \frac{1}{|x|^{n+b-2}}$ as $ x\rightarrow \infty$, we obtain $ u\in L^p(\mathbb{R}^n_+)$, $ \forall p>\frac{n}{n+b-2}$. By $ \eqref{G2}$ we have
\[
| Du(x)| \leq C\int_{\mathbb{R}^n_+} \frac{1}{|x-y|^{n-1}} u(y)^{\tilde{\tau}} dy.
\]
We apply the Hardy-Littlewood-Sobolev inequality to get
\[
\| Du \|_{L^p(\mathbb{R}^n_+)} \leq C \|u^{\tilde{\tau}}\|_{L^{\frac{np}{n+p}}(\mathbb{R}^n_+)}
\]
for any $ p > \frac{n}{n-1}$.
Since $ \tilde{\tau}\frac{np}{n+p} > \frac{n}{n+b-2}$, we obtain $ Du \in L^p(\mathbb{R}^n_+)$. This implies the desired result for all $ p>n$.

\noindent \textbf{Step 2:} \emph{ We prove that $ u \in C^2(\overline{\mathbb{R}^n_+})$.}

We write $ W_R = \{ x=(x',x_n) \ | \ 0 < x_n < R, \ |x'| < R \}$.
For an arbitrarily fixed  $ \bar{x} \in \overline{\mathbb{R}^n_+}$, there is positive constant $ R$, such that $ \bar{x} \in W_{R}$ and $ \text{dist}(\bar{x},\tilde{\p} W_{R}) > 1$, where $ \tilde{\p} W_{R} = \p W_R \backslash \{x \ | \ x_n = 0, \ |x'| < R \}$. We have $ B_{1}(\bar{x})\cap \mathbb{R}^n_+  \subset W_R$.
Take a cutoff function $ \varphi(y) \in C^{\infty}_c(B_{1}(\bar{x}))$ such that $ 0\leq \varphi \leq 1$ in $ B_1(\bar{x})$, $ \varphi =1$ in $ B_{\frac{1}{2}}(\bar{x})$. We define
\[
u_1(x)=\int_{\mathbb R^n_+} K_b(x,y)y^b_n f(u(y))\varphi(y)dy,  \quad x\in B_{\frac{1}{4}}(\bar{x}).
\]
By \cite{Horiuchi95}, we have $ u_1 \in C^{2,\alpha}(\overline{\mathbb R^n_+})$.

Set
\begin{align*}
u_2(x)
& =\int_{\mathbb R^n_+} K_b(x,y)y^b_n f(u(y))(1- \varphi(y)) dy \\
& =\int_{\mathbb R^n_+ \setminus B_{\frac{1}{2}}(\bar{x})} K_b(x,y)y^b_n f(u(y))(1- \varphi(y)) dy, \quad x\in B_{\frac{1}{4}}(\bar{x}).
\end{align*}
Then by the H$ \ddot{o}$lder inequality and $ f(t)=\bar{c}t^{\tilde{\tau}}$, we conclude
\begin{align*}
|D^2 u_2(x)|
& \leq \left \| \frac{1}{|x-y|^n}  \right \|_{L^{q^\ast}(\mathbb R^n_+ \setminus B_{\frac{1}{2}}(\bar{x}))} \| u^{\tilde{\tau}} \|_{L^{p^\ast}(\mathbb R^n_+ \setminus B_{\frac{1}{2}}(\bar{x}))} \\
& \leq C \| u \|^{\tilde{\tau}}_{L^{p}(\mathbb R^n_+ \setminus B_{\frac{1}{2}}(\bar{x}))}
\end{align*}
for $ p > \max \left \{ \frac{n+b+2}{n+b-2}, n \right \} $, where $ p^{\ast}= \frac{p}{\tilde{\tau}}$ and $ \frac{1}{p^{\ast}} + \frac{1}{q^{\ast}} = 1$.
Similarly,
\begin{align*}
|D^3 u_2(x)|
& \leq \left \| \frac{1}{|x-y|^{n+1}} \right \|_{L^{q^\ast}(\mathbb R^n_+ \setminus B_{\frac{1}{2}}(\bar{x}))} \| u^{\tilde{\tau}} \|_{L^{p^\ast}(\mathbb R^n_+ \setminus B_{\frac{1}{2}}(\bar{x}))} \\
& \leq C \| u \|^{\tilde{\tau}}_{L^{p}(\mathbb R^n_+ \setminus B_{\frac{1}{2}}(\bar{x}))}.
\end{align*}
This implies the desired result.\hfill $ \Box$

Since $ u\in C^2(\overline{\mathbb{R}^n_+})$, we may assume that (see \cite{Huang2012})
\begin{equation}\label{s14}
\bar{u}(x',x_n,x_{n+1}) = u\left(x',\sqrt{x^2_n+x^2_{n+1}}\right).
\end{equation}
It follows that,
\begin{equation}\label{sec2}
\left \{
\begin{aligned}
\Delta_{n+1} \bar{u} + \frac{b-1}{x_{n+1}} \bar{u}_{n+1} + \bar{c} \bar{u}(x)^{\tilde{\tau}} &= 0 \quad \text{in}\quad \mathbb R^{n+1}_+ \\
\frac{\p \bar{u}}{\p x_{n+1}} &=0 \quad \text{on} \quad \p \mathbb R^{n+1}_+.
\end{aligned}
\right.
\end{equation}
$ \bar{u}$ is a classical solution to \eqref{sec2}.
\begin{lemma}\label{lem7}
Let $ \bar{u}(x) $ be the positive solution of $ \eqref{sec2} $, and the asymptotic behavior of $ \bar{u}$ at infinity is $ \bar{u}(x) \sim \frac{1}{|x|^{n+b-2}}$. Then $ \bar{u}(x)$ satisfies the corresponding integral equation
\begin{equation}\label{sec3}
\bar{u}(x) = \bar{c} \int_{\mathbb R^{n+1}_+} K_{b-1}(x,y) y^{b-1}_{n+1} \bar{u}(y)^{\tilde{\tau}} dy.
\end{equation}
\end{lemma}
\noindent \emph{Proof.} \ Consider $ R$ large enough and $ |x|>R$, and set $ A_1 := \left\{y\in \mathbb R^{n+1}_+ \ | \ |y-x| \leq \frac{|x|}{2} \right \}$, and $ A_2 := \left\{ y\in \mathbb R^{n+1}_+ \ | \ |y-x| \geq \frac{|x|}{2} \right \}$. Assume that $ \zeta(x) := \bar{c} \int_{\mathbb R^{n+1}_+} K_{b-1}(x,y) y^{b-1}_{n+1} \bar{u}(y)^{\tilde{\tau}} dy$. By the asymptotic behavior of $ \bar{u}$ at infinity, one can easily verify that
\[
\left|\bar{c} \int_{A_1} K_{b-1}(x,y) y^{b-1}_{n+1} \bar{u}(y)^{\tilde{\tau}} dy\right| \leq C\frac{1}{|x|^{n+b}}
\]
and
\begin{align*}
& \left|\bar{c} \int_{A_2} K_{b-1}(x,y) y^{b-1}_{n+1} \bar{u}(y)^{\tilde{\tau}} dy\right| \\
\leq & \ C \int_{A_2 \cap B_{\sqrt{|x|}}} \frac{y^{b-1}_{n+1}}{|x-y^*|^{b-1}} \frac{1}{|x-y|^{n-1}}\left(\frac{1}{1+|y|}\right)^{n+b+2} dy  +  C \int_{A_2 \backslash B_{\sqrt{|x|}}} \frac{1}{|x-y|^{n-1}}\frac{1}{|y|^{n+b+2}} dy \\
\leq & \ \frac{C}{|x|^{n-1+\frac{b-1}{2}}}\int_{A_2 \cap B_{\sqrt{|x|}}} \left(\frac{1}{1+|y|}\right)^{n+b+2} dy + \frac{C}{|x|^{n-1}}\int_{A_2 \backslash B_{\sqrt{|x|}}}\frac{1}{|y|^{n+b+2}} dy \\
\leq & \ \frac{C}{|x|^{n+\frac{b-3}{2}}} + \frac{C}{|x|^{n-1+\frac{b+1}{2}}} \\
\leq & \ \frac{C}{|x|^{n+\frac{b-3}{2}}},
\end{align*}
where $ |x|>R$ and $ B_R = \{x\in \mathbb{R}^{n+1} \ | \ |x|<R \} $. In getting the above inequality, we used $ \forall |y|<\sqrt{|x|}$, $ u(y)\sim \left(\frac{1}{1+|y|}\right)^{n+b-2}$ and
\[
|x-y^*|\geq ||x|-|y^*|| \geq |x|-\sqrt{|x|} \geq \frac{|x|}{2}.
\]
Then we have $ | \zeta(x)| \leq C\frac{1}{|x|^{n+\frac{b-3}{2}}}$. Similarly, we also obtain $ | \nabla \zeta(x)| \leq C\frac{1}{|x|^{n+\frac{b-1}{2}}}$. We know that the asymptotic behavior $ \bar{u}(x) \sim \frac{1}{|x|^{n+b-2}}$ and $ \nabla \bar{u}(x) \sim \frac{1}{|x|^{n+b-1}}$ at $ \infty$. Let $ w=\bar{u} - \zeta$, then $ w$ satisfies the following equation
\[
\Delta w + \frac{b-1}{x_{n+1}} w_{n+1} = 0 \quad \text{in}\quad \mathbb R^{n+1}_+.
\]
Multiplying this identity by $ wx_{n+1}$ and integrating by parts, we get
\begin{align*}
0
& = \lim_{R\rightarrow \infty} \int_{B^+_R(0)} (\Delta w + \frac{b-1}{x_{n+1}} w_{n+1}) wx_{n+1} dx \\
& = - \lim_{R\rightarrow \infty} \int_{B^+_R(0)} x_{n+1}|\nabla w|^2 dx + \lim_{R\rightarrow \infty} \int_{B^+_R(0)} (b-2)ww_{n+1} dx.
\end{align*}
Since the asymptotic properties of $ w$ and $ \nabla w $, the boundary integral term is zero. We deduce that
\begin{align*}
\lim_{R\rightarrow \infty} \int_{B^+_R(0)} x_{n+1}|\nabla w|^2 dx & = \frac{b-2}{2} \lim_{R\rightarrow \infty} \int_{\partial B^+_R(0)} w^2 \nu^{n+1} dS \\
& \leq C \lim_{R\rightarrow \infty} \int_{\partial B^+_R(0)} \frac{1}{|x|^{2n+b-3}} dS \\
& = 0,
\end{align*}
where $ \nu = (\nu^1 , \nu^2 , \cdots , \nu^{n+1}) $ is the outward pointing unit normal vector field. Then we have $ \nabla w = 0$ and $ w = 0$. \hfill $ \Box$

\begin{proposition}\label{pro1}
Let $ f(t)$ satisfy the conditions of Theorem \ref{thm0}. The problem
\begin{equation}\label{pr1}
\left\{
\begin{aligned}
t^b u''(t) + bt^{b-1} u'(t)+t^bf(u(t))&=0 \quad t \in[0,+\infty) \\
u'(0)&= 0
\end{aligned}
\right.
\end{equation}
has no positive solutions of class $ C^2$.
\end{proposition}
In order to prove Proposition \ref{pro1} we make use of the following results. And the idea of the proof from Theorem 3.2 of \cite{EM93}.
\begin{lemma}\label{lem9}
Let $ u(t) \in C^2$ be a positive solution of $ (\ref{pr1})$, then for every $ t>0$ we have
\begin{equation}\label{l3}
tu'(t) + (b-1)u(t) \geq 0.
\end{equation}
\end{lemma}
\noindent $ \emph{Proof.}$ Since $ (\ref{pr1})$ we have
\[
tu'' + bu' + tf(u) =0, \quad t>0,
\]
hence
\[
tu'' + bu'=(tu')' + (b-1)u' \leq 0.
\]
The function $ M(t):= tu'(t)+(b-1)u(t)$ is non-increasing. Now we proceed by contradiction. If there exists $ t_1>0$ such that
\[
M(t_1)=t_1u'(t_1)+(b-1)u(t_1) <0,
\]
then the monotonicity of $ M(t)$ and the positivity of $ u(t)$ imply that
\[
u'(t)\leq u'(t) + \frac{b-1}{t} u(t) \leq \frac{M(t_1)}{t} \ \text{for} \ t\geq t_1.
\]
Integrating the inequality $ u'(t)\leq\frac{M(t_1)}{t}$ on $ (t_1,t)$, we obtain
\[
-u(t_1) \leq u(t)-u(t_1) \leq M(t_1)\ln \left (\frac{t}{t_1}\right).
\]
Letting $ t\rightarrow \infty$ and recalling that $ M(t_1)<0$, we get a contradiction. \hfill $ \Box$

The next lemma contains the necessary a priori estimate for our study.
\begin{lemma}\label{lem10}
Let $ u \in C^2$ be a positive solution of $ (\ref{pr1})$, then
\begin{equation}\label{l7}
u(t) \leq Ct^{-\frac{n+b-2}{2}},
\end{equation}
\begin{equation}\label{l8}
|u'(t)| \leq Ct^{-\frac{n+b}{2}},
\end{equation}
where $ C$ is a positive constant.
\end{lemma}
\noindent $ \emph{Proof.}$ \ We claim $ f(0)=0$. Suppose $ f(0)>0$. Since $ f(t)$ is non-decreasing, we have
\begin{equation}\label{l4}
(t^b u'(t))'=-t^b f(u(t)) \leq -f(0)t^b.
\end{equation}
Integrating $ (\ref{l4})$ on $ (0,t)$ and using the fact that $ u'(0)=0$, we obtain
\begin{equation}\label{l5}
t^b u'(t) \leq -\frac{f(0)}{b+1} t^{b+1}.
\end{equation}
Integrating $ (\ref{l5})$, we get
\[
u(t) \leq -\frac{f(0)}{2(b+1)}t^2 + u(0).
\]
This contradicts $ u(t)>0$ by letting $ t\rightarrow \infty$.

We have $ f(0)=0$. Since $ f(t)$ is non-decreasing, by Remark \ref{r1} there are three cases for $ f$. For the case $ f\equiv 0$, we obtain $ u\equiv 0$ by applying $ (\ref{s3})$. If $ f(t) = 0 \ \text{in} \ [0,t_0]$ for some $ t_0 >0$, then we have $ g(t) = 0 \ \text{in} \ [0,t_0]$ and $ g(t) > 0 \ \text{in} \ (t_0, +\infty)$. This contradicts the non-increasing property of $ g$. Then we consider that $ f(t) > 0$ for every $ t>0$.

Since
\[
(t^b u'(t))' = -t^b f(u(t)) <0, \quad \forall t>0,
\]
$ t^b u'(t)$ is strictly monotonically decreasing. Noticing that $ u'(0)=0$, we see that $ u'(t) \leq 0$, $ \forall t\geq 0$.

Since $ u' \leq 0$ and $ u >0$, we get $ u$ is bounded. Assume that $ u(t)\leq M$ in $ [0,+\infty)$. Since $ g(t)$ is non-increasing, we have
\[
\frac{f(t)}{t^{\tilde{\tau}}} = g(t) \geq g(M) = \frac{f(M)}{M^{\tilde{\tau}}} \quad \forall t\in(0,M].
\]
We can prove $ f(u)\geq Cu^{\tilde{\tau}}$.

Therefore, we have
\begin{equation}\label{l6}
-(t^b u'(t))' = t^b f(u(t)) \geq Ct^b u(t)^{\tilde{\tau}}.
\end{equation}
Integrating $ (\ref{l6})$ on $ (0,t)$, by using $ u'(t)\leq 0$, we have
\[
-t^b u' \geq Ct^{b+1}u(t)^{\tilde{\tau}}.
\]
Integrating $ -u' \geq Ct u(t)^{\tilde{\tau}}$ over $ (0,t)$, we get $ u(t)^{1-\tilde{\tau}} \geq Ct^2 + u(0)^{1-\tilde{\tau}}$. Therefore, we obtain
\[
u(t) \leq Ct^{-\frac{2}{\tilde{\tau}-1}} = Ct^{-\frac{n+b-2}{2}}.
\]
In order to prove $ (\ref{l8})$ it is sufficient to $ (\ref{l3})$ and $ (\ref{l7})$ . \hfill $\Box$ \\

\noindent \textbf{Proof of Proposition \ref{pro1}}

We will proceed by contradiction argument. We define $ F(t)=\int^t_0 f(s) ds$. Multiplying $ (\ref{pr1})$ by $ u$ and integrating by parts on $ (0,t)$, we obtain
\begin{equation}\label{l9}
t^bu(t)u'(t) - \int^t_0 s^b(u'(s))^2 ds  + \int^t_0 s^bf(u(s))u(s) ds = 0.
\end{equation}
On the other hand, multiplying $ (\ref{pr1})$ by $ tu'(t)$ and integrating by parts on $ (0,t)$, we get
\begin{equation}\label{l10}
\frac{1}{2}t^{b+1}(u'(t))^2 + \frac{b-1}{2}\int^t_0 s^b(u'(s))^2 ds  + t^{b+1}F(u(t)) - (b+1)\int^t_0 s^b F(u(s)) ds = 0.
\end{equation}
Using $ (\ref{l7})$ and $ (\ref{l8})$, we conclude that
\[
\lim_{t\rightarrow \infty} t^bu(t)u'(t) = \lim_{t\rightarrow \infty}t^{b+1}(u'(t))^2 =0
\]
and
\[
\int^{\infty}_0 s^b(u'(s))^2 ds < \infty.
\]
Hence by $ (\ref{l9})$, we have
\[
\int^{\infty}_0 s^b(u'(s))^2 ds  = \int^{\infty}_0 s^bf(u(s))u(s) ds < \infty.
\]
We claim that there is a sequence $ t_k\rightarrow \infty$ such that $ t^{b+1}_{k}F(u(t_k))\rightarrow 0$. Suppose not, $ t^{b+1}F(u(t)) \geq C_0 > 0$, $ \forall t > 1$, for some positive constant $ C_0$, we have
\begin{equation}\label{l11}
\frac{C_0}{t} \leq t^{b}F(u(t)) = t^{b}\int^{u(t)}_0 f(s)ds \leq t^bf(u(t))u(t).
\end{equation}
Integrating the inequality $ (\ref{l11})$ on $ (1,+\infty)$, we obtain
\[
\int^{\infty}_1 \frac{C_0}{t} dt \leq \int^{\infty}_1 t^bf(u(t))u(t) dt.
\]
The left-hand side of the inequality is unbounded, but the right-hand side is bounded. This is a contradiction. The claim is proved. Since $ g(t)=\frac{f(t)}{t^{\tilde{\tau}}}$ is non-increasing, we have
\begin{equation}\label{l12}
F(u(t))=\int^{u(t)}_0 \frac{f(s)}{s^{\tilde{\tau}}}s^{\tilde{\tau}}ds \geq \frac{1}{\tilde{\tau}+1}f(u(t))u(t).
\end{equation}
Combining $ (\ref{l9})$, $ (\ref{l10})$ and $ (\ref{l12})$, and taking $ t=t_k$, we get
\begin{align*}
\frac{1}{2}t^{b+1}_k(u'(t_k))^2 + \frac{b-1}{2}\left[t^b_k u(t_k)u'(t_k) + \int^{t_k}_0 s^bf(u(s))u(s) ds\right] \\ + t^{b+1}_k F(u(t_k)) - \frac{(b+1)}{\tilde{\tau}+1}\int^{t_k}_0 s^b f(u(s))u(s) ds \geq 0.
\end{align*}
Letting $ t_k \rightarrow \infty$, we have
\[
\left(\frac{b-1}{2}-\frac{(b+1)}{\tilde{\tau}+1}\right)\int^{\infty}_0 s^b f(u(s))u(s) ds \geq 0.
\]
This contradicts the fact that
\[
\left(\frac{b-1}{2}-\frac{(b+1)}{\tilde{\tau}+1}\right) = \frac{1}{2}\left(\frac{2(b+1)}{n+b}-2\right) < 0.
\]
\hfill $ \Box$\\
\textbf{Proof of Theorem \ref{thm0}}

By Lemma \ref{lem1}, Lemma \ref{lem2} and Lemma \ref{lem3}, we choose the $ x_1 $ direction and prove that $ v $ is symmetric in the $ x_1$ direction. If $ \lam_1 >0$, then $ v$ is symmetric in the direction of $ x_1$. If $ \lam_1 =0 $, then we conclude by continuity that $ v(x) \leq v_0(x)$ for all $ x \in \Sigma_0$. We can also start the moving plane from $ -\infty$ and find a corresponding $ \lam^\prime_1$. If $ \lam^\prime_1 =0$, then we get $ v_0(x) \leq v(x)$ for $ x \in \Sigma_0$. So $ v(x)$ is symmetric with respect to $ T_0$. If $ \lam^\prime_1<0$, an analogue to Lemma \ref{lem3} shows that $v$ is symmetric with respect to $ T_{\lam^\prime_1}$. For the $ x_2,x_3 \cdots x_{n-1}$, we can carry out the procedure as the above. There are two cases for solutions.

\noindent \textbf{Case 1} \  If $ \lam_1 >0 $ or $ \lam^\prime_1 <0$ in some direction for some $ x_0 \in \partial\mathbb{R}^n_+$, we have $ v=v_{\lam_1}$ or $ v=v_{\lam^\prime_1}$. Since $ g$ is non-increasing and by Lemma \ref{lem1}, we get $ g(t) = \bar{c}$. This implies $ f(t)=\bar{c} t^{\frac{n+b+2}{n+b-2}}$. For $ \bar{c}=0$, one has $ u\equiv0$. In the following, we always assume $ \bar{c}>0$. Since $ v$ is regular at the origin, we have the asymptotic behavior of $ u$ at infinity is $ u(x) \sim \frac{1}{|x|^{n+b-2}}$.

According to Lemma \ref{lem8}, we obtain $ u\in C^2(\overline{\mathbb{R}^n_+})$. We know $ u$ satisfies \eqref{s13}. After transformation \eqref{s14}, $ \bar{u}$ satisfies $ n+1$ dimensional equation \eqref{sec2}. Thanks to Lemma \ref{lem7}, we obtain the equivalence between the integral equation \eqref{sec3} and its corresponding differential equation \eqref{sec2}. And we know that $ \bar{u}$ satisfies $ n+1$ dimensional integral equation \eqref{sec3}. Using the above moving plane in integral form for $ \bar{u}$, we obtain $ \bar{u}$ is radially symmetric in the direction of $ x_1,\cdots x_{n-1},x_n $. We define $ (x',x_n)=(x_1,\cdots , x_{n-1},x_n)$. There exists $ p=(p_1,\cdots, p_n) \in \mathbb{R}^{n}$ such that
\[
u(x',|x_n|)= \bar{u}(x',x_n,0)=\bar{u}(\bar{x}',\bar{x}_n,0) = u(\bar{x}',|\bar{x}_n|)
\]
if $\sum^{n}_{i=1}|x_i-p_i|^2 = \sum^n_{i=1}|\bar{x}_i-p_i|^2$. $ u$ is radially symmetric about $ p$. In fact, $ p_n$ must be zero. Otherwise, it follows that
\[
\bar{u}(x',2p_n-x_n,0) = \bar{u}(x',x_n,0) = u(x',|x_n|) = \bar{u}(x',-x_n,0)
\]
It shows that for the fixed $ x'$, $ \bar{u}$ is periodic with respect to $ x_n$ with period $ 2p_n$. Similarly to the proof of Lemma \ref{lem2}, we obtain that $ \bar{u}$ is monotonic with respect to $ p$. Thanks to the decay estimate of $ \bar{u}$ we directly get $ \bar{u} \equiv 0$ and $ u \equiv 0$. We get a contradiction. Thus, we have $ p_n=0$ and $ u$ is radially symmetric about $ p\in\pa\mathbb{R}^{n}_+$. By a result of \cite{ChenLiOu06}, we deduce that $ u(x)=u_{a,p}(x)=\left(\frac{ca}{a^2+|x-p|^2}\right)^{\frac{n+b-2}{2}}$, where $ c=\sqrt{\frac{(n+b)(n+b-2)}{\bar{c}}}$ and $ a\geq0 $.

\noindent \textbf{Case 2} \ Now we suppose that $ \lam_1=\lam^\prime_1=0$ for $  x_1,x_2 \cdots x_{n-1}$ directions and for all $ x_0\in \partial\mathbb{R}^n_+$, then $ v $ and hence $ u $ are radially symmetric in the $  x_1,x_2 \cdots x_{n-1}$ directions. We define $ S_C := \{x\in \mathbb{R}^n_+ \ | \ x_n= C\}$. For any given $ p$, $ q \in S_C$, there exists $ x_0 \in \partial\mathbb{R}^n_+ $ such that $ d(p,x_0) = d(q,x_0)$. Since $ u$ is radially symmetric, we have $ u(p)=u(q)$. According to the arbitrariness of $ p$ and $ q$, the solution $ u$ depends only on $ x_n$. Set $ \tilde{u}(x_n) := u(x)$, and we have:
\begin{equation}
\left\{
\begin{aligned}
x^b_n \tilde{u}^{\prime \prime} + b x^{b-1}_n \tilde{u}^{\prime} +x^b_n f(\tilde{u})&=0 \\
\tilde{u}'(0)&= 0.
\end{aligned}
\right.
\end{equation}
According to Proposition \ref{pro1}, we have $ \tilde{u}\equiv 0$. Hence, $ u(x)\equiv0$. \hfill $ \Box$

\section{ Integral system}
The idea of proving Theorem \ref{thm1} is similar to the proof of Theorem \ref{thm0}. We consider the Kelvin's transform $ w$, $ z$ of $ u$, $ v$ defined by
\[
w(x)= \frac{1}{|x|^{n+b-2}}u \left (\frac{x}{|x|^2} \right ), \ z(x) = \frac{1}{|x|^{n+b-2}}v\left(\frac{x}{|x|^2}\right ).
\]
The definition of $ \Sigma_\lambda$, $ T_\lambda$, $ x^\lambda$ and $ u_\lambda(x)$ in the moving plane method is the same as above.

By a simple calculation, we have
\begin{equation}
\left\{
\begin{aligned}
w(x)&=\int_{\mathbb R^n_+}K_b(x,y)y_n^b f(z(y)|y|^{n+b-2})\frac{1}{|y|^{n+b+2}}dy \quad x \in \mathbb R^n_+ , \\
z(x)&=\int_{\mathbb R^n_+}K_b(x,y)y_n^b g(w(y)|y|^{n+b-2})\frac{1}{|y|^{n+b+2}}dy \quad x \in \mathbb R^n_+ .
\end{aligned}
\right.
\end{equation}
By the definition of  $h(t)=\frac{f(t)}{t^{\tilde{\tau}}}$, $k(t)=\frac{g(t)}{t^{\tilde{\tau}}}$ we deduce that
\begin{equation}
\left\{
\begin{aligned}
w(x)&=\int_{\mathbb R^n_+}K_b(x,y)y_n^b h(z(y)|y|^{n+b-2}) z(y)^{\tilde{\tau}}dy , \\
z(x)&=\int_{\mathbb R^n_+}K_b(x,y)y_n^b k(w(y)|y|^{n+b-2}) w(y)^{\tilde{\tau}}dy .
\end{aligned}
\right.
\end{equation}
Since that $ w$ and $ z$ are continuous and positive in $ \overline{\mathbb R^n_+}$ with possible singularity at the origin. They decay at infinity as $ u(0)|x|^{2-n-b}$ and $ v(0)|x|^{2-n-b}$. We obtain $ w$, $ z\in L^{\frac{2n}{n+b-2}}(\mathbb{R}^n_+ \setminus B_r(0)) \cap L^\infty (\mathbb{R}^n_+ \setminus B_r(0))$ for any $ r>0$.
\begin{lemma}\label{lem4}
\begin{align*}
w(x)-w_\lam(x)=& \int_{\Sigma_\lam} (K_b(x,y)-K_b(x^\lam,y)) y_n^b \times \\
&[h(z(y)|y|^{n+b-2})z(y)^{\tilde{\tau}}- h(z(y^\lam) |y^{\lam}|^{n+b-2}) z(y^\lam)^{\tilde{\tau}} ] dy,
\end{align*}
\begin{align*}
z(x)-z_\lam(x)=& \int_{\Sigma_\lam} (K_b(x,y)-K_b(x^\lam,y)) y_n^b \times \\
&[k(w(y)|y|^{n+b-2})w(y)^{\tilde{\tau}}- k(w(y^\lam) |y^{\lam}|^{n+b-2}) w(y^\lam)^{\tilde{\tau}} ] dy.
\end{align*}
\end{lemma}
\noindent \emph{Proof.} The proof is similar to that of
Lemma \ref{lem1} and is omitted. \hfill $ \Box $
\begin{lemma}\label{lem5}
Under the conditions of Theorem \ref{thm1}, there exists $ \lambda_0 > 0$ such that for all $ \lambda \geq \lambda_0$, we have $ w_\lambda(x) \geq w(x)$ and $ z_\lambda(x) \geq z(x)$  for all $ x \in \Sigma_\lambda$.
\end{lemma}
\noindent \emph{Proof.} \ The proof is similar to Lemma \ref{lem2}. We denote by $ \Sigma^w_\lam = \{y\in \Sigma_\lam \ | \ w(y) > w_\lam(y)\}$ and $ \Sigma^z_\lam = \{y\in \Sigma_\lam \ | \ z(y) > z_\lam(y)\}$. By the Hardy-Littlewood-Sobolev inequality, we get
\begin{align*}
\|w-w_\lam\|_{L^q(\Sigma^w_\lam)}
&\leq C\left(\int_{\Sigma_\lam^z}z^{\frac{2n}{n+b-2}}(y)dy\right)^{\frac{2}{n}}\|z-z_\lam\|_{L^q(\Sigma^z_\lam)},\\
\|z-z_\lam\|_{L^q(\Sigma^z_\lam)}
&\leq C\left(\int_{\Sigma_\lam^w}w^{\frac{2n}{n+b-2}}(y)dy\right)^{\frac{2}{n}}\|w-w_\lam\|_{L^q(\Sigma^w_\lam)},
\end{align*}
where $ q>\frac{n}{n-2}$. The above two inequalities imply
\begin{align}
&\|w-w_\lam\|_{L^q(\Sigma^w_\lam)} \nonumber \\
&\leq C\left(\int_{\Sigma_\lam^w}w^{\frac{2n}{n+b-2}}(y)dy\right)^{\frac{2}{n}}\left(\int_{\Sigma_\lam^z}z^{\frac{2n}{n+b-2}}(y)dy\right)^{\frac{2}{n}}\|w-w_\lam\|_{L^q(\Sigma^w_\lam)}, \label{s15}
\end{align}
and
\begin{align}
&\|z-z_\lam\|_{L^q(\Sigma^z_\lam)} \nonumber \\
&\leq
C\left(\int_{\Sigma_\lam^w}w^{\frac{2n}{n+b-2}}(y)dy\right)^{\frac{2}{n}}\left(\int_{\Sigma_\lam^z}z^{\frac{2n}{n+b-2}}(y)dy\right)^{\frac{2}{n}}\|z-z_\lam\|_{L^q(\Sigma^z_\lam)}. \label{s16}
\end{align}
Since $ w(x)$, $ z(x)\in L^{\frac{2n}{n+b-2}}(\mathbb{R}^n_+ \setminus B_r(0))$ for any $ r>0$, we choose $ \lam_0$ large enough such that
\[
C\left(\int_{\Sigma_\lam}w^{\frac{2n}{n+b-2}}(y)dy\right)^{\frac{2}{n}}\left(\int_{\Sigma_\lam}z^{\frac{2n}{n+b-2}}(y)dy\right)^{\frac{2}{n}} \leq \frac{1}{2}
\]
for all $ \lam \geq \lam_0$. Then we have
\[
\|w-w_\lam\|_{L^q(\Sigma^w_\lam)} = \|z-z_\lam\|_{L^q(\Sigma^z_\lam)} = 0
\]
for all $ \lam \geq \lam_0$. So we get the desired result.  \hfill $ \Box$

We define $ \lam_1 = \inf \{ \lam \ | \ w(x)\leq w_\mu(x) \ \text{and} \ z(x)\leq z_\mu(x), \ \forall \mu \geq \lam, \  x\in\Sigma_\mu \}$.
\begin{lemma}\label{lem6}
If $ \lam_1>0$, then $ w(x)\equiv w_{\lam_1}(x)$ and $ z(x)\equiv z_{\lam_1}(x)$ for all $ x\in \Sigma_{\lam_1}$.
\end{lemma}
\noindent \emph{Proof.} \ Suppose that $  z(x)\not\equiv z_{\lam_1}(x)$, then we can infer from Lemma \ref{lem4} that $ w< w_{\lam_1}$ in the interior of $ \Sigma_{\lam_1}$ and this further implies $ z< z_{\lam_1}$ in the same area. If one can show that for $ \eps$ sufficiently small so that $ \forall \lam \in (\lam_1-\eps,\lam_1]$, there holds
\begin{equation}\label{s17}
C\left(\int_{\Sigma_\lam^w}w^{\frac{2n}{n+b-2}}(y)dy\right)^{\frac{2}{n}}\left(\int_{\Sigma_\lam^z}z^{\frac{2n}{n+b-2}}(y)dy\right)^{\frac{2}{n}} \leq \hlf,
\end{equation}
then by \eqref{s15} and \eqref{s16}, we have
\[
\|w-w_\lam\|_{L^q(\Sigma^w_\lam)} = \|z-z_\lam\|_{L^q(\Sigma^z_\lam)} = 0,
\]
and therefore $ \Si^w_\lam$ and $ \Si^z_\lam$ must be measure zero. Since $ w$ and $ z$ are continuous, we deduce that  $ \Si^w_\lam$ and $ \Si^z_\lam$ are empty. This contradicts the definition of $ \lam_1$.

Now we verify inequality \eqref{s17}. Since $ w$, $ z \in L^{\frac{2n}{n+b-2}}(\mathbb{R}^n_+ \setminus B_r(0))$, for any small $ \eta >0$, we can choose $ R$ sufficiently large so that
\[
C\left(\int_{\mathbb{R}^n_+ \backslash B_R}w^{\frac{2n}{n+b-2}}(y)dy\right)^{\frac{2}{n}} \leq \eta, \qquad C\left(\int_{\mathbb{R}^n_+ \backslash B_R}z^{\frac{2n}{n+b-2}}(y)dy\right)^{\frac{2}{n}} \leq \eta.
\]
We fix this $ R$ and then show that the measure of $ \Si^w_\lam \cap B_R$ and $ \Si^z_\lam \cap B_R$ are sufficiently small for $ \lam$ close to $ \lam_1$. The rest part of the proof is similar to that of Lemma \ref{lem3} and is omitted. \hfill $\Box$ \\
\textbf{Proof of Theorem \ref{thm1}}

By Lemma \ref{lem4}, Lemma \ref{lem5} and Lemma \ref{lem6}, we choose the $ x_1 $ direction and prove that $ w$, $ z$ are symmetric in the $ x_1$ direction. If $ \lam_1 >0$, then $ w$, $ z$ are symmetric in the direction of $ x_1$. If $ \lam_1 =0 $, then we conclude by continuity that $ w(x) \leq w_0(x)$, $ z(x) \leq z_0(x)$ for all $ x \in \Sigma_0$. We can also start the moving plane from $ -\infty$ and find a corresponding $ \lam^\prime_1$. If $ \lam^\prime_1 =0$, then we get $ w_0(x) \leq w(x)$, $ z_0(x) \leq z(x)$ for $ x \in \Sigma_0$. So $ w(x)$ and $ z(x)$ are symmetric with respect to $ T_0$. If $ \lam^\prime_1<0$, an analogue to Lemma \ref{lem6} shows that $ w$ and $ z$ are symmetric with respect to $ T_{\lam^\prime_1}$. For the $ x_2,x_3 \cdots x_{n-1}$, we can carry out the procedure as the above. There are two cases for solutions.

\noindent \textbf{Case 1} \  If $ \lam_1 >0 $ or $ \lam^\prime_1 <0$ in some direction for some $ x_0 \in \partial\mathbb{R}^n_+$, we have $ w=w_{\lam_1}$, $ z=z_{\lam_1}$ or $ w=w_{\lam^\prime_1}$, $ z=z_{\lam^\prime_1}$. Since $ h$ and $ k$ are non-increasing and by the Lemma \ref{lem4}, then we get $ h(t) = c_1$ and $ k(t)=c_2$, where $ c_1$ and $ c_2$ are positive constants. We have $ f(t)=c_1 t^{\frac{n+b+2}{n+b-2}}$ and $ g(t)=c_2 t^{\frac{n+b+2}{n+b-2}}$. Since $ w$ and $ z$ are regular at the origin, we have the asymptotic behavior of $ u$ and $ v$ at infinity are $ u \sim \frac{1}{|x|^{n+b-2}}$, $ v \sim \frac{1}{|x|^{n+b-2}}$.

Similarly, we can prove $ u$, $ v \in C^2(\overline{\mathbb{R}^n_+})$. Thus, we take the transformation
\[
\bar{u}(x',x_n,x_{n+1}) = u\left(x',\sqrt{x^2_n+x^2_{n+1}}\right) \quad \text{and}   \quad \bar{v}(x',x_n,x_{n+1}) = v\left(x',\sqrt{x^2_n+x^2_{n+1}}\right).
\]
It is easy to verify that $ \bar{u}$ and $ \bar{v}$ satisfy equation
\begin{equation}\label{s20}
\left \{
\begin{aligned}
\Delta_{n+1} \bar{u} + \frac{b-1}{x_{n+1}} \bar{u}_{n+1} + c_1 \bar{v}(x)^{\tilde{\tau}} & = 0 \quad \text{in}\quad \mathbb R^{n+1}_+ \\
\frac{\p \bar{u}}{\p x_{n+1}} & = 0 \quad \text{on} \quad  \p \mathbb R^{n+1}_+.
\end{aligned}
\right.
\end{equation}
\begin{equation}\label{s21}
\left \{
\begin{aligned}
\Delta_{n+1} \bar{v} + \frac{b-1}{x_{n+1}} \bar{v}_{n+1} + c_2 \bar{u}(x)^{\tilde{\tau}} &= 0 \quad \text{in}\quad \mathbb R^{n+1}_+ \\
\frac{\p \bar{v}}{\p x_{n+1}} & = 0 \quad \text{on} \quad  \p \mathbb R^{n+1}_+.
\end{aligned}
\right.
\end{equation}
Exactly as the proof of Lemma \ref{lem7}, using the asymptotic behavior of $ \bar{u}$ and $ \bar{v}$ at infinity, one can easily deduce that
\begin{equation}\label{s18}
\bar{u}(x) = c_1 \int_{\mathbb R^{n+1}_+} K_{b-1}(x,y) y^{b-1}_{n+1} \bar{v}(y)^{\tilde{\tau}} dy,
\end{equation}
\begin{equation}\label{s19}
\bar{v}(x) = c_2 \int_{\mathbb R^{n+1}_+} K_{b-1}(x,y) y^{b-1}_{n+1} \bar{u}(y)^{\tilde{\tau}} dy.
\end{equation}
We establish the equivalence between the integral equation \eqref{s18}, \eqref{s19} and its corresponding differential equation \eqref{s20}, \eqref{s21}. Using the above moving plane method for the $ n+1$ dimensional integral equations \eqref{s18} and \eqref{s19}, we obtain $ \bar{u}$ and $ \bar{v}$ are radially symmetric in the direction of $ x_1,\cdots x_{n-1},x_n $. There exists $ p \in \pa\mathbb{R}^{n}_+$ such that
\[
u(x',|x_n|)= \bar{u}(x',x_n,0)=\bar{u}(\bar{x}',\bar{x}_n,0) = u(\bar{x}',|\bar{x}_n|)
\]
and
\[
v(x',|x_n|)= \bar{v}(x',x_n,0)=\bar{v}(\bar{x}',\bar{x}_n,0) = v(\bar{x}',|\bar{x}_n|).
\]
if $\sum^{n}_{i=1}|x_i-p_{i}|^2 = \sum^n_{i=1}|\bar{x}_i-p_{i}|^2$. Thus, $ u$ and $ v$ are radially symmetric about $ p$. By the results of \cite{ChenLiOu05} and \cite{ChenLiOu06} we obtain $ u(x)=\left(\frac{ca}{a^2+|x-p|^2}\right)^{\frac{n+b-2}{2}}$ and $ v(x)=\left(\frac{\tilde{c}a}{a^2+|x-p|^2}\right)^{\frac{n+b-2}{2}}$ where $ p \in \pa\mathbb{R}^{n}_+$.

\noindent \textbf{Case 2} \ (The idea of the proof from Theorem 3.2 of \cite{EM93}) Now we suppose that $ \lam_1=\lam^\prime_1=0$ for $  x_1,x_2 \cdots x_{n-1}$ directions and for all $ x_0\in \partial\mathbb{R}^n_+$, then $ w $, $ z$ and hence $ u$, $ v$ are radially symmetric in the $  x_1,x_2 \cdots x_{n-1}$ directions. We define $ S_C := \{x\in \mathbb{R}^n_+ \ | \ x_n= C\}$. Similarly, we obtain the solutions $ u$ and $ v$ depend only on $ x_n$. Set $ \tilde{u}(x_n) := u(x)$, $ \tilde{v}(x_n) := v(x)$ and we have:
\begin{equation}\label{o1}
\left\{
\begin{aligned}
x^b_n \tilde{u}^{\prime \prime} + b x^{b-1}_n \tilde{u}^{\prime} +x^b_n f(\tilde{v}) &=0 \\ \tilde{u}'(0)&=0,
\end{aligned}
\right.
\end{equation}
\begin{equation}\label{o2}
\left\{
\begin{aligned}
x^b_n \tilde{v}^{\prime \prime} + b x^{b-1}_n \tilde{v}^{\prime} +x^b_n g(\tilde{u}) &=0 \\ \tilde{v}'(0)&=0.
\end{aligned}
\right.
\end{equation}
Exactly as the proof of Lemma \ref{lem9}, using $ f(t)\geq 0$ and $ g(t)\geq 0$ we deduce that
\begin{equation}\label{o3}
t \tilde{u}'(t) + (b-1)\tilde{u}(t) \geq 0,
\end{equation}
\begin{equation}\label{o4}
t\tilde{v}'(t) + (b-1)\tilde{v}(t) \geq 0.
\end{equation}
Similarly to the proof of Lemma \ref{lem10}, we obtain
\[
(b-1) \tilde{u}(t) \geq -t\tilde{u}'(t) \geq C t^2 \tilde{v}^{\tilde{\tau}}(t),
\]
\[
(b-1) \tilde{v}(t) \geq -t\tilde{v}'(t) \geq C t^2 \tilde{u}^{\tilde{\tau}}(t).
\]
Solving these inequalities, we get for all $ t>0$,
\begin{equation}\label{o5}
\tilde{u}(t) \leq Ct^{-\frac{n+b-2}{2}}, \qquad \tilde{v}(t) \leq Ct^{-\frac{n+b-2}{2}}.
\end{equation}
Since \eqref{o3} and \eqref{o4}, we have
\begin{equation}\label{o6}
|\tilde{u}'(t)| \leq Ct^{-\frac{n+b}{2}}, \qquad |\tilde{v}'(t)| \leq Ct^{-\frac{n+b}{2}}.
\end{equation}
Multiplying \eqref{o1} by $ \tilde{v}$ and \eqref{o2} by $ \tilde{u}$ and integrating by parts on $ (0,t)$, we get
\begin{equation}\label{o7}
t^b\tilde{u}'(t)\tilde{v}(t) - \int^t_0 x^b_n\tilde{u}'\tilde{v}' dx_n = - \int^t_0 x^b_n f(\tilde{v})\tilde{v}dx_n,
\end{equation}
\begin{equation}\label{o8}
t^b\tilde{u}(t)\tilde{v}'(t) - \int^t_0 x^b_n\tilde{u}'\tilde{v}' dx_n = - \int^t_0 x^b_n g(\tilde{u})\tilde{u}dx_n.
\end{equation}
Using the fact that \eqref{o5} and \eqref{o6}, we deduce that
\begin{equation}
\lim_{t\rightarrow \infty} t^b\tilde{u}'(t)\tilde{v}(t) = \lim_{t\rightarrow \infty} t^b\tilde{u}(t)\tilde{v}'(t) =0,
\end{equation}
\[
\int^{\infty}_0 x^b_n\tilde{u}'\tilde{v}' dx_n < \infty.
\]
Hence by \eqref{o7} and \eqref{o8}, we have
\begin{equation}\label{o9}
\int^{\infty}_0 x^b_n\tilde{u}'\tilde{v}' dx_n = \int^{\infty}_0 x^b_n f(\tilde{v})\tilde{v} dx_n =\int^{\infty}_0 x^b_n g(\tilde{u})\tilde{u}dx_n.
\end{equation}
We define $ F(t)=\int^t_0 f(s)ds$ and $ G(t)=\int^t_0 g(s)ds$. By multiplying  \eqref{o1} by $ x_n\tilde{v}'$ and \eqref{o2} by $ x_n\tilde{u}'$ and integrating by parts on $ (0,t)$, we obtain
\[
t^{b+1}\tilde{u}'(t)\tilde{v}'(t) - \int^t_0 x^b_n\tilde{u}'\tilde{v}' dx_n - \int^t_0 x^{b+1}_n \tilde{u}'\tilde{v}'' dx_n = -t^{b+1}F(\tilde{v}) + (b+1)\int^t_0 x^b_n F(\tilde{v})dx_n,
\]
\[
\int^t_0 x^{b+1}_n \tilde{u}'\tilde{v}'' dx_n +  b \int^t_0 x^b_n\tilde{u}'\tilde{v}' dx_n = -t^{b+1}G(\tilde{u}) + (b+1)\int^t_0 x^b_n G(\tilde{u})dx_n.
\]
Hence, we get
\[
t^{b+1}\tilde{u}'(t)\tilde{v}'(t) + (b-1) \int^t_0 x^b_n\tilde{u}'\tilde{v}' dx_n + t^{b+1}(F(\tilde{v}) + G(\tilde{u})) - (b+1)\int^t_0 x^b_n (F(\tilde{v}) + G(\tilde{u})) dx_n = 0 .
\]
As in the proof of Proposition \ref{pro1}, there is a sequence $ t_k \rightarrow \infty$ such that
\begin{equation}\label{l13}
t^{b+1}_k F(\tilde{v}(t_k)) \rightarrow 0, \qquad  t^{b+1}_k G(\tilde{u}(t_k)) \rightarrow 0.
\end{equation}
Using the fact that \eqref{o5} and \eqref{o6}, we have
\begin{equation}\label{l15}
\lim_{t\rightarrow \infty} t^{b+1}\tilde{u}'(t)\tilde{v}'(t) =0.
\end{equation}
And it is easy to see that
\begin{equation}\label{l14}
F(\tilde{v}(t)) \geq \frac{1}{\tilde{\tau}+1}f(\tilde{v}(t))\tilde{v}(t), \qquad G(\tilde{u}(t)) \geq \frac{1}{\tilde{\tau}+1}g(\tilde{u}(t))\tilde{u}(t).
\end{equation}
By taking $ t=t_k$, we have
\[
t^{b+1}_k\tilde{u}'(t_k)\tilde{v}'(t_k) + (b-1) \int^{t_k}_0 x^b_n\tilde{u}'\tilde{v}' dx_n + t^{b+1}_k (F(\tilde{v}) + G(\tilde{u})) - (b+1)\int^{t_k}_0 x^b_n (F(\tilde{v}) + G(\tilde{u})) dx_n = 0 .
\]
Letting $ t_k \rightarrow \infty$ and using \eqref{o9}, \eqref{l13}, \eqref{l15} and \eqref{l14}, we get
\[
\left(\frac{b-1}{2}-\frac{b+1}{\tilde{\tau}+1}\right) \int^{\infty}_0 (x^b_n f(\tilde{v})\tilde{v} + x^b_n g(\tilde{u})\tilde{u}) dx_n \geq 0.
\]
Since $ \frac{b-1}{2} - \frac{b+1}{\tilde{\tau}+1} < 0$, we get a contradiction. We obtain $ u(x)\equiv0$ and $ v(x)\equiv0$. \hfill $ \Box$\\

\noindent \textbf{Acknowledgments.} \ The work of the author is sponsored by Shanghai Rising-Star Program 19QA1400900. The author would like to thank Professor Huang Genggeng for his valuable suggestions. The author was supported by National Natural Science Foundation of China under Grant 11871160.

\bibliographystyle{plain}
\bibliography{citeliouville}

{\em Addresse and E-mail:}
\medskip

{\em Yating Niu}

{\em School of Mathematical Sciences}

{\em Fudan University}

{\em ytniu19@fudan.edu.cn}

\end{document}